\documentclass[12pt]{amsart}
\usepackage{amscd}

\def\NZQ{\Bbb}               

\def\ZZ{{\NZQ Z}}

\def\frk{\frak}               

\def\mm{{\frk m}}

\def\opn#1#2{\def#1{\operatorname{#2}}} 
\opn\chara{char}
\opn\length{\ell}
\opn\pd{pd}
\opn\rk{rk}
\opn\projdim{proj\,dim}
\opn\rank{rank}
\opn\depth{depth}
\opn\grade{grade}
\opn\height{height}
\opn\embdim{emb\,dim}
\opn\codim{codim}

\opn\Tr{Tr}
\opn\bigrank{big\,rank}
\opn\superheight{superheight}\opn\lcm{lcm}
\opn\trdeg{tr\,deg}%
\opn\reg{reg}
\opn\lreg{lreg}
\opn\div{div}
\opn\Div{Div}
\opn\cl{cl}
\opn\Cl{Cl}
\opn\Spec{Spec}
\opn\Supp{Supp}
\opn\supp{supp}
\opn\Sing{Sing}
\opn\Ass{Ass}
\opn\Ann{Ann}
\opn\Rad{Rad}
\opn\Soc{Soc}
\opn\Ker{Ker}
\opn\Coker{Coker}
\opn\Im{Im}
\opn\Hom{Hom}
\opn\Tor{Tor}
\opn\Ext{Ext}
\opn\End{End}
\opn\Aut{Aut}
\opn\id{id}

\opn\nat{nat}
\opn\pff{pf}
\opn\Pf{Pf}
\opn\GL{GL}
\opn\SL{SL}
\opn\mod{mod}
\opn\ord{ord}
\opn\aff{aff}
\opn\con{conv}
\opn\relint{relint}
\opn\st{st}
\opn\lk{lk}
\opn\cn{cn}
\opn\core{core}
\opn\vol{vol}
\opn\link{link}
\opn\star{star}
\opn\gr{gr}

\def\poly#1#2#3{#1[#2_1,\dots,#2_{#3}]}
\def\pot#1#2{#1[\kern-0.28ex[#2]\kern-0.28ex]}

\opn\dirlim{\underrightarrow{\lim}}
\opn\inivlim{\underleftarrow{\lim}}
\let\dirsum=\oplus

\let\iso=\cong

\let\Dirsum=\bigoplus

\let\mcone= * 
\let\to=\rightarrow
\def\namedTo#1{\overset{#1}{\longrightarrow}}
\let\To=\longrightarrow
\def\Implies{\ifmmode\Longrightarrow \else
     \unskip${}\Longrightarrow{}$\ignorespaces\fi}
\def\implies{\ifmmode\Rightarrow \else
     \unskip${}\Rightarrow{}$\ignorespaces\fi}
\def\iff{\ifmmode\Longleftrightarrow \else
     \unskip${}\Longleftrightarrow{}$\ignorespaces\fi}

\let\:=\colon
\newtheorem{Theorem}{Theorem}[section]
\newtheorem{Lemma}[Theorem]{Lemma}
\newtheorem{Corollary}[Theorem]{Corollary}
\newtheorem{Proposition}[Theorem]{Proposition}
\newtheorem{Remark}[Theorem]{Remark}

\newtheorem{Example}[Theorem]{Example}

\let\epsilon\varepsilon
\let\kappa=\varkappa
\opn\inii{in}
\opn\inim{inm}
\opn\set{set}
\def\pnt{{\raise0.5mm\hbox{\large\bf.}}}

\begin{document}

\title{Single Spot Ideals of Codimension~3 and 
long Bourbaki Sequences}
\author{Yukihide Takayama}
\address{Yukihide Takayama, Department of Mathematical
Sciences, Ritsumeikan University, 
1-1-1 Nojihigashi, Kusatsu, Shiga 525-8577, Japan}
\email{takayama@se.ritsumei.ac.jp}


\def\Coh#1#2{H_{\mm}^{#1}(#2)}
\def\cal{\bold}  

\newcommand{\AppTh}{Theorem~\ref{approxtheorem} }
\def\da{\downarrow}
\newcommand{\ua}{\uparrow}
\newcommand{\namedto}[1]{\buildrel\mbox{$#1$}\over\rightarrow}
\newcommand{\bdel}{\bar\partial}
\newcommand{\proj}{{\rm proj.}}

\begin{abstract}
Let $K$ be a field and $S = K[x_1,\ldots, x_n]$ be a polynomial ring.
A single spot ideal $I \subset S$ is a graded
ideal whose local cohomology $\Coh{i}{S/I}$, $i<\dim S/I$ and 
$\mm = (x_1,\ldots, x_n)$, only has
non-trivial value $N$, a finite length module, at $i = \depth S/I$. We
consider characterization of single spot ideals in terms of (long)
Bourbaki sequences. The codimension~2 case has been fairly well
investigated.  In this paper, we focus on the codimension~3 case.
\end{abstract}

\maketitle
\section*{Introduction}
Let $S = \poly{K}{x}{n}$ be a polynomial ring over a field $K$ with 
the standard grading and let
$\mm = (x_1,\ldots, x_n)$. All the modules and ideals in this paper
are graded. 
A finitely generated $S$-module $M$ is
called a {\em generalized Cohen-Macaulay $(CM)$ module} if the local
cohomoloy module $\Coh{i}{M}$ has finite length for all $i<\dim(M)$. 
A ring $R$ is a {\em generalized CM ring} if it is a generalized CM
$R$-module. An ideal $I\subset S$ is called a {\em generalized CM
ideal} if $S/I$ is a generalized CM ring. If a generalized CM module
$M$ satisfies $\dim M = \dim S$, it is called {\em maximal}.

In this paper, we are interested in generalized CM ideals. 
In particular, single spot ideals. An ideal $I\subset S$ is called a 
{\em single spot ideal} of type $(t,N)$ where $t = \depth S/I$ and $N$ 
is a finite length
$S$-module if the local cohomology only has a non-trivial value $N$ at
dimension $t$, i.e.,
\[
    \Coh{i}{S/I} = \left\{
			\begin{array}{ll}
				0 & \mbox{for all $i<\dim S/I$ with $i\neq t$} \\
                                N & \mbox{if $i=t$}
			\end{array}
                   \right.
\]
Let $I\subset S$ be a generalized CM ideal of $\codim I = r$ $(r\geq 2)$. 
Then by Corollary~1.3 \cite{HT2}
we have a long Bourbaki sequence
\[
	0\To F_{r-1} \To \cdots \To F_1 \To M \To I \To 0
\]
with $S$-free modules $F_i$ $(i=1,\ldots, r-1)$ and $M$ is a maximal generalized CM 
module whose local cohomology is as follows 
\begin{equation}
\label{localcohomologies}
	\Coh{i}{M}  \iso 
\left\{
\begin{array}{ll}
  \Coh{i}{I} \iso \Coh{i-1}{S/I}   & \mbox{if }i < n-r+1 \\
0                              & \mbox{if } i=n-r+1\\
\end{array}
\right.
\end{equation}
In this sence, the ideal $I$ is approximated by $M$.
Notice that the value of $\Coh{i}{M}$ for $i=n-r+2,\ldots, n-1$ are
irrelevant to this approximation. However,
the construction given in the proof of Corollary~1.3
\cite{HT2} (and also Lemma~1.3 \cite{Ama1} in a slightly different
situation) always makes the module
$M$ such that $\Coh{i}{M} =0$ for $i=n-r+2,\ldots, n-1$.
In this paper, we are interested in
long Bourbaki sequences with approximation modules $M$ such that
$\Coh{i}{M}$ $(i=n-r+2,\ldots, n-1)$  are not always trivial,
and study the case of codimension~3, namely the case of 
$\Coh{n-1}{M}=N$ where $N$ is a non-trivial finite length module.
Notice that in the case of codimension~2 we always have $\Coh{n-1}{M}=0$.

First of all, we will give a characterization of a maxmal generalized
CM module $M$ whose local cohomology is
\begin{equation}
\label{lc-condition-codim3}
\Coh{i}{M} = \left\{
		\begin{array}{ll}
		    K  & \mbox{if $i=t+1$}\\
                    N  & \mbox{if $i=n-1$} \\
                    0  & \mbox{if $i<n-2$, $i\ne t+1$}
		\end{array}
             \right.
\end{equation}
in terms of the first syzygy of $M$. See Theorem~\ref{theorem:main1}.
Then we consider the special case of $M = E_{t+1}\dirsum E_{n-1}(d)$ where
$E_j$ denotes the $j$th sygyzy module of the field $K$ over $S$. 
We will use the notation $M(d)$, for a graded module $M$ and $d\in\ZZ$,
such that $M(d)_i= M_{d+i}$ (the $d+i$th component of $M$) for all $i\in\ZZ$.
Our question is 
how we can construct a long Bourbaki sequence 
\begin{equation}
\label{the:bourbaki-sequence}
   0 \To F\To G \To E_{t+1}\dirsum E_{n-1}(d) \To I(c) \To 0
\end{equation}
of {\em non-trivial} type. Here a {\em trival type} construction is as follows. First
construct a long Bourbaki sequence
\[
   0 \To F' \namedTo{f} G' \namedTo{f} E_{t+1} \namedTo{\phi} I(c) \To 0
\]
according to the method given in the proof of Corollary~1.3 \cite{HT2}
(or Lemma~1.3 \cite{Ama1}). Then make the direct sum
\[
   0 \To F'\dirsum K_n(d) \namedTo{f\dirsum\partial_n}
 G'\dirsum K_{n-1}(d) \namedTo{g\dirsum\partial_{n-1}}
 E_{t+1}\dirsum E_{n-1}(d)  \namedTo{\phi} I(c) \To 0
\]
where $(K_\bullet, \partial_\bullet)$ is the Koszul complex of the sequence 
$x_1,\ldots, x_n$ over $S$. We will denote a base $x_{i_1}\wedge\cdots\wedge x_{i_k}$
($1\leq i_1 < \cdots < i_k\leq n$) of the Koszul complex of the sequence
$x_1,\ldots, x_n$ by $e_{\{i_1,\ldots,i_k\}}$ or $e_{i_1\cdots i_k}$.   
Notice that, in the trivial type Bourbaki sequence
$E_{n-1}(d)$ does not contribute to $I$ via $\phi$.

It is well known that a (short) Bourbaki sequence $0 \to F \to M \To I
\to 0$ is constructed by finding '(graded) basic elements' in $M$.
See \cite{Bour,Kunz} for the standard basic element theory
and \cite{HT2,F} for graded version. 
However, there is no comparative notion  for long Bourbaki 
sequences. We give a simple answer to this problem in the case of 
$E_{t+1}\dirsum E_{n-1}(d)$ (and $M= E_{t+1}$).
We will give a characterization of long Bourbaki sequences
(\ref{the:bourbaki-sequence}) in terms of elemets from $K_{t+1}\dirsum
K_{n-1}$ (from $K_{t+1}$) 
and from $K_{n-t-1}\dirsum K_{1}$ (from $K_{n-t-1}$) satisfying certain
conditions, which suggests a consruction of the long Bourbaki sequences.  See
Theorem~\ref{theorem:construction} and~\ref{theorem:construction2}. 
In particular, non-trivial type
construction is characterized by an additional condition on the
elements from $K_{t+1}\dirsum K_{n-1}$
(Theorem~\ref{theorem:non-trivialcond}).

However, the existence of a long Bourbaki sequence
(\ref{the:bourbaki-sequence}) only means that $I$ is a single spot
ideal of codimension less than or equal to $3$. We give a numerial
condition to assure $\codim I = 3$. See
Theorem~\ref{theorem:numerical}.  Finally, we give some examples. 

For a module $M$ ($\ne K$), 
we will denote the $j$th syzygy module over $S$ by $\Omega_j(M)$.
Also we use two kinds of duals,
$(-)^\vee = \Hom_S(-,K)$ and $(-)^*= \Hom_S(-, S(-n))$.

\section{Approximation Modules of Single Spot Ideals of type $(t, K(-c))$}
In this section, we consider 
approximation modules $M$ of codimension~3 single spot ideals $I\subset S$
of type $(t, K(-c))$ in long Bourbaki sequences
\[
   0\To F \To G \To  M \To I(c) \To 0.
\]
If we restrict ourself to the case of $\Coh{n-1}{M}=0$, we have 
$M = E_{t+1}\dirsum H$ for some free $S$-module $H$ 
according to Herzog, Takayama \cite{HT2} and Amasaki \cite{Ama1}.
We will now consider the general case.

\begin{Theorem}
\label{theorem:main1}
Let $M$ be a maximal generalized CM module over $S$ and 
consider its first sygyzy:
\[
   0 \To \Omega_{1}(M) \To F \To M \To 0
\]
and let $g_1,\ldots, g_l$ be a minimal set of generators 
of $\Omega_{1}(M)$ whose
degrees are $a_1,\ldots, a_l$.
Also let $N$ be a finite length module over $S$,
which may be $0$.
Then the following are equivalent.
\begin{enumerate}
 \item [$(i)$]
   For $t \leq n-4$, we have
\[ \Coh{i}{M} =
	\left\{
	\begin{array}{ll}
	   K   & \mbox{if } i=t+1 \\
	0 & \mbox{if } i\leq n-2, i\ne t+1\\
	   N   & \mbox{if } i=n-1 
	\end{array}
      \right.
\]
\item [$(ii)$]
\begin{enumerate}
\item [$(a)$] $\Omega_1(M) \iso \Dirsum_{i=1}^lS(-a_i)/E_{t+3}$, and
\item [$(b)$] $\Omega_1(M)^* \iso N^\vee + F^*/M^*$ 
\end{enumerate}
\end{enumerate}
\end{Theorem}
%
%
\begin{proof}
We first prove $(i)$ to $(ii)$. Let $F_\bullet$ be a minimal free resolution of $M$
over $S$:
\[
 F_\bullet : 
	0\To F_{n-t-1} \namedTo{\varphi_{n-t-1}}
	     F_{n-t-2} \namedTo{\varphi_{n-t-2}}
	         \cdots\namedTo{\varphi_3}
             F_2       \namedTo{\varphi_2}
	     F_1       \namedTo{\varphi_1} 
             F_0       \namedTo{\varphi_0}
               M       \To 0.
\]
By taking the dual, we have 
\[
 0 \To F_0^*      \namedTo{\varphi_1^*} 
       F_1^*      \namedTo{\varphi_2^*}
       F_2^*      \namedTo{\varphi_3^*}
       \cdots     \namedTo{\varphi_{n-t-2}^*} 
       F_{n-t-2}^*\namedTo{\varphi_{n-t-1}^*}
       F_{n-t-1}^*\To 0.
\]
Then by local duality the  $j$th cohomology of this complex is
\[
   \Ext_S^j(M, S(-n)) \iso \Coh{n-j}{M}^\vee
   =   \left\{
	\begin{array}{ll}
                N^\vee & \mbox{if $j=1$} \\
                0      & \mbox{if $j\geq 2, j\ne n-t-1$}\\
		K      & \mbox{if $j=n-t-1$} 
	\end{array}
\right.
\] 
for $j\geq 1$. Thus 
\[
0\To \Im\varphi_2^* \To F_2^* \To\cdots\To F_{n-t-1}^*\To K \To 0
\]
is exact and $F_2^*$ to $F_{n-t-1}^*$ part is a begining of a minimal
free resolution of $K$, which is isomorphic to a begining of the Koszul
complex $(K_\bullet,\partial_\bullet)$ of the sequence $x_1,\ldots,
x_n$. Namely,
\[
     F_{n-t-1}^* \iso K_0,\quad\ldots\quad, F_2^* \iso K_{n-t-3} \qquad\mbox{and}\quad 
     \Im\varphi_2^* \iso E_{n-t-2}.
\]
On the other hand, we have
$N^\vee \iso \Ker\varphi_2^*/\Im\varphi_1^*$ 
and 
$E_{n-t-2}\iso \Im\varphi_2^* \iso  F_1^*/\Ker\varphi_2^*$.
Now set $U := \Coker\varphi_1^* = F_1^*/\Im\varphi_1^*$. Then
\[
U/N^\vee \iso (F_1^*/\Im\varphi_1^*)/(\Ker\varphi_2^*/\Im\varphi_1^*)
 \iso  F_1^*/\Ker\varphi_2^* = E_{n-t-2}.
\]
Thus  we have
\[
      0 \To N^\vee \To U \To E_{n-t-2} \To 0
\]
Taking the dual, we have
\[
    0\To E_{n-t-2}^* \To U^* \To (N^\vee)^*.
\]
Since $N$ has finite length, $N^\vee$ has also finite length by Matlis
duality, so that $(N^\vee)^* =0$. Also $E_{n-t-2}^* \iso E_{t+3}$
by selfduality of Koszul complex. Thus we have $U^* \iso E_{t+3}$. 
Then by dualizing  the exact sequence
\begin{equation}
\label{prop:new1-seq1}
  0\To M^* \To F_0^* \namedTo{\varphi_1^*}F_1^* \To U \To 0
\end{equation}
we have
\[
	0\To E_{t+3} \To F_1 \namedTo{\varphi_1} F_0 \To M\To 0.
\]
This proves $(ii) (a)$. 
Now from the short exact sequence
\[
    0 \To \Omega_1(M) \To F \namedTo{\varphi} M \To 0
\]
we have the long exact sequence 
\[
   0 \To M^* \namedTo{\varphi^*}
         F^* \To 
         \Omega_1(M)^* \To
         N^\vee \To 0
\]
since we have $\Ext_S^1(M,S(-n))\iso \Coh{n-1}{M}^n = N^\vee$ by local
duality. This proves $(ii) (b)$.

Next we prove $(ii)$ to $(i)$. By $(ii)(a)$ we have a $S$-free resolution
of $M$:
\[
   0\To K_n\namedTo{\partial_n}\cdots\namedTo{\partial_{t+4}}
          K_{t+3}\namedTo{\partial_{t+3}}
        F_1\namedTo{\varphi_1} F_0 \namedTo{\varphi_0} M \To 0
\]
where $F_0$ and $F_1$ are $S$-free modules. By taking the dual, we have
the complex
\[
    0\To 
      M^*    \namedTo{\varphi_0^*}
      F_0^* \namedTo{\varphi_1^*} F_1^* \namedTo{\partial_{t+3}^*}
	K_{t+3}^* \namedTo{\partial_{t+4}^*}\cdots\namedTo{\partial_n^*}
       K_n^* \To 0
\]
Then by local duality and selfduality of Koszul
complex we compute 
\[
    \Coh{i}{M} \iso \Ext_S^{n-i}(M,S(-n))^\vee
                =\left\{
			\begin{array}{ll}
			    K   & \mbox{if $i = t+1$} \\
                            0   & \mbox{if $i\leq n-2$, $i\ne t+1$}
			\end{array}
                \right.
\]
Now by dualizing the exact sequence
\[
    K_{t+3}\namedTo{\partial_{t+3}} F_1 \namedTo{\varphi_1} \Ker\varphi_0 \To 0
\]
we have 
\[
   0\To (\Ker\varphi_0)^* \To F_1^* \namedTo{\partial_{t+3}^*} K_{t+3}^*,
\]
so that we have $\Omega_1(M)^* =(\Ker\varphi_0)^*\iso \Ker(\partial_{t+3}^*)$.
Then by the condition $(ii)(b)$ we compute
\begin{eqnarray*}
\Coh{n-1}{M}^\vee 
      &\iso &\Ext_S^{1}(M, S(-n)) = \Ker \partial_{t+3}^* / \Im \varphi_1^*
\\
      &\iso & \Omega_1(M)^* / (F_0^*/\varphi_0^*(M^*))
\\
      &  = & N^\vee
\end{eqnarray*}
as required.
\end{proof}

\begin{Corollary}
Let $M$ be a maximal generalized CM module satisfying Theorem~\ref{theorem:main1} (i).
Then its minimal free resolution is in the form of 
\[
0\to K_n\namedTo{\partial_n}
     K_{n-1}\namedTo{\partial_{n-1}}\cdots\namedTo{\partial_{t+4}}
     K_{t+3}\namedTo{\partial_{t+3}}
     F_1\namedTo{\varphi}
     F_0 \to M \to 0.
\]
\end{Corollary}
Notice that $F_1$ is a $S$-free module containing a submodule isomorphic to
$E_{t+3}$.

\begin{Example}\label{example:1}
Let $M = E_{t+1}\dirsum E_{n-1}$. 
Then $\Coh{i}{M}$ is as in Theorem~\ref{theorem:main1} $(i)$
with $N=K$.
Since $M$ has a minimal free resolution
\[
  0\To K_n \namedTo{\partial_n}\cdots\namedTo{\partial_{t+4}}
       K_{t+3}\namedTo{\partial_{t+3}}
       K_{t+2}\dirsum K_{n} \namedTo{\partial_{t+2}\dirsum \partial_n}
       K_{t+1}\dirsum K_{n-1}\namedTo{\partial_{t+1}\dirsum \partial_{n-1}} M \To 0,
\]
we have $\Omega_1(M) = E_{t+2}\dirsum E_{n} \iso G/E_{t+3}$
where $G = K_{n+2} \dirsum E_{n}$. Thus we have Theorem~\ref{theorem:main1} $(ii) (a)$.
On the other hand, 
we have $\Omega_1(M)^* = E_{t+2}^*\dirsum E_{n}^* \iso E_{n-t-1}\dirsum S$
by selfduality of Koszul complex and $E_{n} \iso S(-n)$.
Again by selfduality we have
$(K_{t+1}\dirsum K_{n-1})^*/M^*
=
(K_{t+1}\dirsum K_{n-1})^*/(E_{t+1}\dirsum E_{n-1})^*
= (K_{n-t-1}\dirsum K_1)/ (E_{n-t}\dirsum E_2)
= (K_{n-t-1}/E_{n-t})\dirsum (K_1/E_2)
= E_{n-t-1} \dirsum E_1
= E_{n-t-1}\dirsum \mm$.
Thus $\Omega_1(M)^*/((K_{t+1}\dirsum K_{n-1})^*/M^*) \iso S/\mm \iso K \iso K^\vee$ and 
we obtain Theorem~\ref{theorem:main1} $(ii)(b)$.
\end{Example}

\section{Long Bourbaki sequences with Approximation
Module $E_{t+1}\dirsum E_{n-1}(d)$}
In the last chapter, we considered approximation modules $M$ satisfying 
the condition of Theorem~\ref{theorem:main1}(i). We now focus on a 
special case of $M = E_{t+1}\dirsum E_{n-1}(d)$, and investigate the 
long Bourbaki sequences.

We will use the following well known result frequently without
refering it. First of all, we will give a proof for the readers'
convenience.

\begin{Lemma} 
\label{lemma:wellknown}
For the $t$-th syzygy module of $K$ over $S$, we 
have 
\[
\rank E_t =\binom{n-1}{t-1}
\]
\end{Lemma}
\begin{proof}
Let $K_i$ be the $i$th Koszul complex. Then we have
an exact sequence
\[
  0\To K_n\namedTo{\partial_n} K_{n-1}\namedTo{\partial_{n-1}}\cdots\namedTo{\partial_t} K_t \To E_t \To 0 
\]
so that 
\[
	\rank(E_t) =\sum_{i=t}^{n}(-1)^{i-t}\rank(K_i)
		   = (-1)^t\sum_{i=t}^{n}(-1)^i\binom{n}{i} 
= (-1)^{t-1}\sum_{i=0}^{t-1}(-1)^i\binom{n}{i}.
\]
Now set $\alpha(n,t) := \rank(E_t)$. 
By a straightforward calculation, we have
$\alpha(n,t) - \alpha(n-1, t) = \alpha(n-1,t-1)$.
Thus we have 
\begin{eqnarray*}
\alpha(n,t) &= &\alpha(n-1,t) + \alpha(n-1, t-1) \\
            & =& (-1)^{t-1}\sum_{i=0}^{t-1}(-1)^i\binom{n-1}{i}
              + (-1)^{t-2}\sum_{i=0}^{t-2}(-1)^i\binom{n-1}{i} 
	     = \binom{n-1}{t-1}.
\end{eqnarray*}
as required.
\end{proof}

\subsection{Characterization of long Bourbaki Sequences}

For $J, K \subset [n]=\{1,2,\ldots, n\}$ with $J\cap K =\emptyset$ we
define $\sigma(J, K) = (-1)^i$ where $i = \sharp\{ (j,k)\in J\times K
\mid j > k\}$. Then we have $x_J\wedge x_K = \sigma(J,K)x_{J\cup K}$.

Now a long Bourbaki sequence with approximation module $E_{t+1}\dirsum
E_{n-1}(d)$ is characterized by suitable sequences from
$K_{t+1}\dirsum K_{n-1}$ and its dual. Namely,

\begin{Theorem}
\label{theorem:construction}
Following are equivalent.
\begin{enumerate}
\item [$(i)$] We have a long Bourbaki sequence
\[
 0 \To \Dirsum_{i=1}^p S(-a_i) 
   \To \Dirsum_{i=1}^q S(-b_i) 
   \To  E_{t+1}\dirsum E_{n-1}(d)
   \To I(c) \To 0 \quad (exact)
\]
where  $I \subset S$ is a graded ideal.
\item [$(ii)$]
We have 
$\beta_1,\ldots,\beta_q \in K_{t+1}\dirsum K_{n-1}(d)\backslash E_{t+2}\dirsum E_n(d)$
and $\varphi = (a,b)$ with
$a \in {\cal A}$ and $b\in {\cal B}$,  
where 
\begin{eqnarray*}
{\cal A}& =& \left\langle
		\sum_{j=1}^{n-t}
			(-1)^{j+1}\sigma(L\backslash\{i_j\}, 
([n]\backslash L)\cup\{i_j\})
				x_{i_j} e^*_{([n]\backslash L)\cup\{i_j\}}
		\mid L = \{i_1,\ldots,i_{n-t}\}\subset [n]
            \right\rangle \\
{\cal B} &= &\langle
		(-1)^ix_j e^*_{[n]\backslash\{i\}} -
		(-1)^jx_i e^*_{[n]\backslash\{j\}}
		\mid 1\leq i < j \leq n
            \rangle,
\end{eqnarray*}
such that 
\begin{enumerate}
\item [$(a)$] 
      $\varphi :  K_{t+1}\dirsum K_{n-1}(d)\to S(-n)$ is 
      a degree '$n+c$' homomorphism and 
      $\Ker(\varphi)
	= \langle\beta_1,\ldots\beta_q\rangle + E_{t+2}\dirsum E_n(d)$, and 
\item [$(b)$] we we have the following diagram, with $p = q - n + 2 -\binom{n-1}{t}$
\[
\begin{CD}
     @.    0              @.            0         \\
 @.        @VVV                         @VVV      \\
 0  @>>>  \Ker Res(\beta)         @>>>  \Ker\beta  @>>> 0\\
 @.        @VVV                         @VVV             \\
 0  @>>>  \Dirsum_{i=1}^pS(-a_i)  @>>>  \Dirsum_{i=1}^qS(-b_i) \\
 @.       @V{Res(\beta)}VV              @V{\beta}VV \\
0 @>>> \langle\beta_1,\ldots,\beta_q\rangle
      \cap E_{t+2}\dirsum E_n(d) 
                               @>>>  \langle\beta_1,\ldots,\beta_q\rangle \\
 @.        @VVV                         @VVV \\
     @.    0            @.               0   \\
\end{CD}
\]
where $\beta(g_i) = \beta_i$ for all $i$ with $g_1,\ldots, g_q$ the free basis 
of $\Dirsum_{i=1}^q S(-b_i)$, and $Res(-)$ denotes the restriction of maps.
\end{enumerate}
\end{enumerate}
In this case, we have $I = \varphi(K_{t+1}\dirsum K_{n-1}(d))(-c)$
\end{Theorem}
\begin{proof}
We first prove $(ii)$ to $(i)$. 
First notice that,  by the selfduality of Koszul complex, we have
\begin{eqnarray*}
E_i & \iso &  E^*_{n-i+1}\iso\partial^*_{n-i+1}(E^*_{n-i+1}) \\
& =   & \left\{ \sum_{k=1}^{i}(-1)^{k+1} 
                   \sigma(J\backslash \{j_k\}, [n]-(J\backslash\{j_k\}))
                    x_{j_k} e^*_{[n]-(J\backslash\{j_k\})}
: J = \{j_1,\ldots, j_i\}\subset [n]
\right\}.\\
\end{eqnarray*}
See \cite{BH} Chapter 1.6. Thus ${\cal A} = \partial^*_{t+1}(E_{t+1}^*)$ and 
${\cal B} = \partial^*_{n-1}(E_{n-1}^*)$. Then, there exists 
$\bar{a} \in E_{t+1}^*$ and $\bar{b}\in E_{n-1}^*$ such that 
$a = \bar{a}\circ \partial_{t+1}$ and $b = \bar{b}\circ\partial_{n-1}$.
Then by $(a)$ we have the diagram
\begin{equation}
\label{construction:proof-1}
\begin{CD}
0 @>>> \langle\beta_1,\ldots,\beta_q\rangle
       + E_{t+2}\dirsum E_n(d) @>>> K_{t+1}\dirsum K_{n-1}(d) @>{\varphi}>> S(c) \\
@.     @V{\bdel}VV                 @V{\bdel}VV                         @|   \\
0 @>>> \langle\bdel\beta_1,\ldots,\bdel\beta_q\rangle
                              @>>> E_{t+1}\dirsum E_{n-1}(d)@>{\phi}>>  S(c) \\
@.      @VVV                        @VVV                                @.  \\
   @.    0                    @.     0                       @.            
\end{CD}
\end{equation}
where $\phi = (\bar{a},\bar{b})$ and $\bdel := \partial_{t+1}\dirsum \partial_{n-1}(d)$.
On the other hand, we have 
\begin{eqnarray*}
 \Ker\bdel\circ\beta & = & 
   \{ \sum_{i=1}^q h_i g_i \mid \bdel(\sum_{i=1}^q h_i\beta_i) = 0, h_i\in S \} 
\\                   & = &
   \{ \sum_{i=1}^q h_i g_i \mid \sum_{i=1}^q h_i\beta_i\in\Ker\bdel = E_{t+2}\dirsum E_n(d), h_i \in S\}
\\                   & = & 
   \{ \sum_{i=1}^q h_i g_i \mid \beta(\sum_{i=1}^q h_ig_i)
                                 \in\langle\beta_1,\ldots, \beta_q\rangle
				\cap  E_{t+2}\dirsum E_n(d)\}
\\                   & \iso & \Dirsum_{i=1}^p S(-a_i)
\end{eqnarray*}
where the last isomorphism is by $(b)$. 
Notice that let $u\in\Dirsum_{i=1}^q S(-b_i)$ be such that
$\beta(u)\in\langle\beta_1,\ldots,\beta_q\rangle\cap
E_{t+2}\dirsum E_{n}(d)$. 
Then $u$ must be in $\Dirsum_{i=1}^p S(-a_i)$. In fact,
by (b) we can choose 
$v\in \Dirsum_{i=1}^pS(-a_i)$ such that
$Res(\beta)(v) = \beta(u)$. Thus $u- v \in\Ker\beta 
\iso \Ker Res(\beta)$, and we have $u \in v + \Ker Res(\beta)
\subset \Dirsum_{i=1}^pS(-a_i)$ as required.
Then by (\ref{construction:proof-1}) we obtain
\[
  0 \To \Dirsum_{i=1}^p S(-a_i) \To \Dirsum_{i=1}^q S(-b_i)
	\namedTo{\bdel\circ\beta} E_{t+1}\dirsum E_{n-1}(d) 
        \namedTo{\phi} S(c)\qquad (exact)
\]
and since $p = q - n + 2 -\binom{n-1}{t}$ we know that $\rank\;\Im\psi =1$
so that we have 
\[
  0 \To \Dirsum_{i=1}^p S(-a_i) \To \Dirsum_{i=1}^q S(-b_i)
	\namedTo{\bdel\circ\beta} E_{t+1}\dirsum E_{n-1} 
        \namedTo{\phi} I(c)\To 0\qquad (exact)
\]
for the ideal $I := \phi(E_{t+1}\dirsum E_{n-1}(d))(-c)
                =  \varphi(K_{t+1}\dirsum K_{n-1}(d))(-c)$
as required.

Next we prove $(i)$ to $(ii)$. Given a long Bourbaki sequence
\[
  0 \To \Dirsum_{i=1}^p S(-a_i) 
 	\namedTo{f} \Dirsum_{i=1}^q S(-b_i)
	\namedTo{g}E_{t+1}\dirsum E_{n-1}(d) 
        \namedTo{\phi} I(c)\To 0\qquad (exact)
\]
with a graded ideal $I\subset S$. Then we have $p = q - n + 2 + \binom{n-1}{t}$
since $1 = \rank I(c) = \rank E_{t+1}\dirsum E_{n-1}(d) - q + p$. Also since  $I(c) \subset S(c)$ we have 
\begin{eqnarray*}
\phi & \in &  \Hom_S(E_{t+1}\dirsum E_{n-1}(d), S(c)) \\
& =  &  \Hom_S(E_{t+1}\dirsum E_{n-1}(d), S(-n))(n+c)\\
& =  & (E_{t+1}^*\dirsum E_{n-1}^*(-d))(n+c) 
\iso(\partial_{t+1}^*(E_{t+1}^*)\dirsum \partial_{n-1}^*(E_{n-1}^*(-d)))(n+c).\\
\end{eqnarray*}
Thus there exists a unique $(a,b) \in 
\partial_{t+1}^*(E_{t+1}^*)\dirsum \partial_{n-1}^*(E_{n-1}^*(-d))$ such that 
$\phi\circ\bdel = (a,b)$. Now we set $\varphi = (a,b)$.
Then $E_{t+2}\dirsum E_{n}(d) + N = \bdel^{-1}(\Ker\phi) (\subset K_{t+1}\dirsum
K_{n-1}(d))$ for some module $N(\ne 0)$.
Let $\{\beta_1,\ldots, \beta_q\}$ be a minimal set of generators of $N$. 
Then we have the following diagram:
\[
\begin{CD}
   @.                                           @.      0        @. \\
@.       @.                                             @AAA     @. \\
   @.    \Dirsum_{i=1}^qS(-b_i)                 @>{g}>> E_{t+1}\dirsum E_{n-1}(d)
@>{\phi}>>          I(c) @>>> 0 \\
@.       @V{\beta}VV                                  @A{\bdel}AA
                    @|   \\
0 @>>> \langle\beta_1,\ldots,\beta_q\rangle @>>> K_{t+1}\dirsum K_{n-1}(d) @>{\varphi=(a,b)}>> I(c) \\
@.       @VVV                                         @.               @. \\
    @.     0                                    @.           @.           \\
\end{CD}
\]
Then we have obtained $\beta_1,\ldots, \beta_q \in K_{t+1}\dirsum K_{n-1}(d)$ and $(a,b)\in
{\cal A}\times {\cal B}$ satisfying the condition $(a)$ and the numerical condition on $p$, $q$, $t$ and  $n$ in $(b)$. Also we have 
$I = \varphi(K_{t+1}\dirsum K_{n-1}(d))(-c)$.

Now since $\Ker g =\Ker (Res(\bdel)\circ\beta) = \Im f =
\Dirsum_{i=1}^p S(-a_i)$ we readily have the following diagram:
\[
\begin{CD}
   @.    0               @.       0          @.    \\
@.      @VVV                      @VVV           @.\\
0  @>>> \Ker\beta\circ f @>>>     \Ker\beta @>>> 0 \\
@.      @VVV                      @VVV           @.\\
0 @>>> \Dirsum_{i=1}^pS(-a_i) 
                         @>{f}>> \Dirsum_{i=1}^qS(-b_i) 
@>{g}>> 
E_{t+1}\dirsum E_{n-1}(d)\\
@.       @.                       @V{\beta}VV     @| \\
0 @>>> \langle\beta_1,\ldots,\beta_q\rangle\cap E_{t+2}\dirsum E_n(d)
                         @>>>     \langle\beta_1,\ldots,\beta_q\rangle
                                                  @>{Res(\bdel)}>>
                                                  E_{t+1}\dirsum E_{n-1}(d)\\
@.       @.                       @VVV             @. \\
    @.                   @.       0               @. \\
\end{CD}
\]
Notice that since $\Ker\beta \subset\Ker (Res(\bdel)\circ\beta) 
= \Ker g = \Im f$ we have the exactness of the first row.
Since 
$\Im (\beta\circ f) = \beta (\Ker(Red(\bdel)\circ \beta)) =  \Ker Res(\bdel) 
= \langle\beta_1,\ldots,\beta_q\rangle\cap E_{t+2}\dirsum E_n(d)$, 
we have a well-defined surjection
\[
   \beta\circ f :  \Dirsum_{i=1}^p S(-a_i)
 \To \langle\beta_1,\ldots,\beta_q\rangle
					\cap E_{t+2}\dirsum E_n(d).
\]
Thus we obtained the diagram of $(b)$ as required.
\end{proof}

For any codimension~3 single spot ideal of type $(t, K(-c))$, there exists
a long Bourbaki sequence with approximation module $E_{t+1}\dirsum H$ where 
$H$ is a $S$-free module \cite{HT2,Ama1}. The following is the case 
of $H=0$, which can be proved with the same idea as that of  
Theorem~\ref{theorem:construction}.
\begin{Theorem}
\label{theorem:construction2}
Following are equivalent.
\begin{enumerate}
\item [$(i)$] We have a long Bourbaki sequence
\[
 0 \To \Dirsum_{i=1}^p S(-a_i) 
   \To \Dirsum_{i=1}^q S(-b_i) 
   \To  E_{t+1}
   \To I(c) \To 0 \quad (exact)
\]
where  $I\subset S$ is a graded ideal.
\item [$(ii)$]
We have 
$\beta_1,\ldots,\beta_q \in K_{t+1}\backslash E_{t+2}$
and $\varphi\in {\cal A}$
where 
\[
{\cal A} = \left\langle
		\sum_{j=1}^{n-t}
			(-1)^{j+1}\sigma(L\backslash\{i_j\},
 ([n]\backslash L)\cup\{i_j\})
				x_{i_j} e^*_{([n]\backslash L)\cup\{i_j\}}
		\mid L=\{i_1,\ldots,i_{n-t}\} \subset [n] \right\rangle
\]
such that 
\begin{enumerate}
\item [$(a)$] 
       $\varphi : K_{t+1} \to S(-n)$ defines a degree '$n+c$' homomorphism and 
	$\Ker\varphi = \langle\beta_1,\ldots\beta_q\rangle + E_{t+1}$, and 
\item [$(b)$] we we have the following diagram, with $p = q + 1 -\binom{n-1}{t}$
\[
\begin{CD}
    @.    0              @.    0        @.      \\
@.       @VVV                  @VVV           @.\\
0 @>>> \Ker Res(\beta) @>>> \Ker\beta @>>>  0 \\
@.       @VVV                  @VVV           @.\\
0 @>>> \Dirsum_{i=1}^pS(-a_i)
                         @>>>  \Dirsum_{i=1}^qS(-b_i)
                                        @.       \\
@.       @V{Res(\beta)}VV       @V{\beta}VV    @. \\
0 @>>> \langle\beta_1,\ldots,\beta_q\rangle\cap E_{t+2}
                         @>>>  \langle\beta_1,\ldots,\beta_q\rangle
                                               @. \\
@.        @VVV                  @VVV           @. \\
     @.    0             @.     0        @.        \\
\end{CD}
\]
where $\beta(g_i) = \beta_i$ for all $i$ with $g_1,\ldots, g_q$ the free basis 
of $\Dirsum_{i=1}^q S(-b_i)$.
\end{enumerate}
In this case, we have $I = \varphi(K_{t+1})(n-c)$.
\end{enumerate}
\end{Theorem}

\begin{Corollary}
\label{cor:depthzero}
There is no codimension~3 single spot ideal $I$ of type $(0, K(-c))$ with 
approximation module $E_{t+1}$ if $t=0$. 
\end{Corollary}
\begin{proof}
Assume that there exists a codimension~3 single spot ideal
$I \subset S$ of type $(0, K(-c))$
fitting into a long Bourbaki sequence
\[
   0 \To F \To G \To E_1 \To I(c) \To 0.
\]
Then by Theorem~\ref{theorem:construction2}, there exist
$\beta_1,\ldots, \beta_q \in K_1\backslash E_2$ and 
$(0\ne) \varphi \in {\cal A}$
such that $\langle\beta_1,\ldots,\beta_q\rangle + E_2 = \Ker(\varphi : K_1\to S(-n))$. Since ${\cal A} =  \langle x_1e^*_1 +\cdots +  x_ne^*_n\rangle$
we must have $\beta_1,\ldots, \beta_q\in E_2$, 
a contradiction.
\end{proof}

\begin{Remark}
In fact, the above Corollary holds for any codimension $\geq 2$. 
Here is an outline of the proof. We consider a straightforward extension of 
Theorem~\ref{theorem:construction2}(a) to any codimension $r$.
In this case we consider long Bourbaki sequences
\[
   0\To F_{r-1} \To F_{r-2}\To\cdots\To F_{1}\To E_1 \namedTo{\varphi} I(c)\To 0.
\]
and $\varphi$ is determined by a nonzero element from
${\cal A} = \langle x_1e^*_1 +\cdots +  x_ne^*_n\rangle$.
\end{Remark}

\subsection{Long Bourbaki sequences of non-trivial type}

Let $n\geq 4$ and $t\leq n-4$ and 
consider a long Bourbaki sequence
\begin{equation}
\label{bourbakiseq0}
    0 \To F\namedTo{f} G \namedTo{g} E_{t+1}\dirsum E_{n-1}(d) \namedTo{\phi} I(c) \To 0
\end{equation}
with $c\in\ZZ$ and  $S$-free modules $F$ and $G$.
From this sequence, we construct the following
diagram where the second row is the minimal free resolution of $M =
E_{t+1}\dirsum E_{n-1}(d)$ and the third row is the mapping cone
$C(\alpha,\beta)$ of a chain map $\alpha\dirsum\beta$, which is a free
resolution of $I(c)$. 
\begin{equation}
\label{mappingcone}
\begin{array}{ccccccccl}
    &     &     &           &            &         &             &      0        &        \\
    &     &     &           &            &         &             &    \da        &        \\
    &  0  & \To &  F        &\namedTo{f} & G       & \namedTo{g} & \Ker\phi      & \to 0  \\
    &     &     & \alpha\da &            &\beta\da &             &  \da          & \      \\
0\to K_n\namedTo{\partial_n}\cdots\namedTo{\partial_{t+4}}
    &K_{t+3} &\namedTo{\partial_{t+3}}&K_{t+2}
& \namedTo{\begin{array}{c}
\partial_{t+2}\\ 
\dirsum\partial_n
\end{array}}
 & K_{t+1} & \namedTo{\begin{array}{c} \partial_{t+1}\\
\dirsum\partial_{n-1}
\end{array}} & M & \to 0 \\
    &       &    &\dirsum K_n(d)&    &\dirsum K_{n-1}(d)&    &\hspace*{0.3cm}\da\phi  &  \\
0\to K_n\namedTo{\partial_n}\cdots\namedTo{\partial_{t+4}} 
& K_{t+3} &\namedTo{\zeta} & K_{t+2}&\namedTo{\rho} & K_{t+1}&  \namedTo{\begin{array}{c}
\phi\circ(\partial_{t+1}\\ \dirsum\partial_{n-1})
\end{array}} &I(c) &\to 0 \\
  &\dirsum F&    &\dirsum K_n(d)&    &\dirsum K_{n-1}(d)&    &  \da  &  \\
  &         &    &\dirsum G  &    &               &    &   0   &  \\
\end{array}
\end{equation}
where
\[
  \begin{array}{cccc}
     \rho:  & K_{t+2}\dirsum K_n(d) \dirsum G & \To  & K_{t+1}\dirsum K_{n-1}(d) \\
      & (a, b, c)  & \longmapsto & (\partial_{t+2}(a),\partial_n(b))- \beta(c)  
  \end{array}
\]
and 
\[
  \begin{array}{cccc}
     \zeta:  & K_{t+3} \dirsum F & \To  & K_{t+2}\dirsum K_{n}(d)\dirsum G \\
      & (a, b)  & \longmapsto & (\partial_{t+3}(a), 0, f(b)) - \alpha(b)  
  \end{array}
\]
Let $p_1 : K_{t+1}\dirsum K_{n-1}(d) \to K_{t+1}$ and $p_2 : K_{t+1}\dirsum
K_{n-1}(d)\to K_{n-1}(d)$ be the first and the second projections. 
From the diagram (\ref{mappingcone})  we know 
    $\Ker\phi = \Im g = (\partial_{t+1}\dirsum \partial_{n-1})\circ \beta(G)$
and then by considering the ranks of the modules in the short exact sequence
\begin{equation}
\label{kerVarphi}
   0 \To \Ker\phi
     \To  E_{t+1}\dirsum E_{n-1}(d)
     \namedTo{\phi} I(c)
     \To 0
\end{equation}
we have
\begin{equation}
\label{rankcond1} 
    \rank(\Ker\varphi) = n-2 +\binom{n-1}{t}.
\end{equation}
On the other hand, we have 
\begin{equation}
\label{inclusion}
(\partial_{t+1}\circ p_1\circ\beta)(G)
         \dirsum (\partial_{n-1}\circ p_2 \circ \beta)(G)
  \supset
(\partial_{t+1}\dirsum\partial_{n-1})\beta(G)
 = \Ker\phi.
\end{equation}
Thus we have
\begin{equation}
\label{rankcond2}
 \rank {\cal I}_{t+1} +\rank{\cal I}_{n-1} \geq n-2 + \binom{n-1}{t}.
\end{equation}
where  ${\cal I}_{t+1} := (\partial_{t+1}\circ p_1\circ \beta)(G) (\subseteq E_{t+1})$ and 
${\cal I}_{n-1} := (\partial_{n-1}\circ p_2\circ \beta)(G) (\subseteq E_{n-1}(d))$.
Since $\rank E_{t+1} = \binom{n-1}{t}$ and $\rank E_{n-1}(d) = n-1$ we know from 
(\ref{rankcond2}) that 
\[
(\rank {\cal I}_{t+1}, \rank {\cal I}_{n-1}) 
= (\binom{n-1}{t} - 1, n-1), 
  (\binom{n-1}{t}, n-1), \mbox{ or }
  (\binom{n-1}{t}, n-2).
\]
Under this situation, we have
\begin{Lemma}
\label{chainmapBETA}
Following are equivalent.
\begin{enumerate}
\item [$(i)$]  long Bourbaki sequence $(\ref{bourbakiseq0})$ is 
               of non-trivial type
\item [$(ii)$] For any free basis $\{m_i\}_{i=1}^q$ of $G$
there exists an index $i$ such that
$\beta(m_i) \notin K_{t+1}$ and $\beta(m_i)\notin K_{n-1}(d)$.
\end{enumerate}
\end{Lemma}
\begin{proof}
We will prove $(i)$ to $(ii)$. 
We assume that for all $i$ we have either
$\beta(m_i)\in K_{t+1}$ or $\beta(m_i)\in K_{n-1}(d)$ and 
will deduce  a contradiction.
First of all, we have equality in (\ref{inclusion}), and then
from (\ref{rankcond2}) we have 
\[
    n-2 + \binom{n-1}{t} = \rank {\cal I}_{t+1} +  \rank {\cal I}_{n-1}.
\]
Thus, we have 
$(\rank {\cal I}_{t+1}, \rank {\cal I}_{n-1}) = 
(\binom{n-1}{t}-1, n-1)$ or $(\binom{n-1}{t}, n-2)$.
Also, since
$\Ker\phi = {\cal I}_{t+1}\dirsum {\cal I}_{n-1}$, we have by (\ref{kerVarphi})
\begin{equation}
\label{prop:I(c)}
I(c) \iso (E_{t+1}/{\cal I}_{t+1})\dirsum (E_{n-1}(d)/{\cal I}_{n-1})
\end{equation}
\begin{description}
\item [case $(\rank {\cal I}_{t+1}, \rank {\cal I}_{n-1}) = (\binom{n-1}{t}-1, n-1)$:]
Since we have
$\rank E_{n-1}(d)/{\cal I}_{n-1} =\rank E_{n-1}(d) - \rank {\cal I}_{n-1} = 0$,
$E_{n-1}/{\cal I}_{n-1}$ is 0 or a torsion-module. But since $I(c)$ is torsion
free, we must have $E_{n-1}(d) = {\cal I}_{n-1}$ by (\ref{prop:I(c)}). Thus 
$\Ker\phi = {\cal I}_{t+1}\dirsum E_{n-1}(d)$ and then the Bourbaki sequence 
(\ref{bourbakiseq0}) must be of trivial-type
\[
     0\To F'\dirsum K_n \namedTo{f'\dirsum\partial_n}
          G'\dirsum K_{n-1}\namedTo{g'\dirsum\partial_{n-1}}
          E_{t+1}\dirsum E_{n-1}(d)\namedTo{\phi} I(c)\To 0
\]
where $0\to F' \to G' \to {\cal I}_{t+1} \to 0$ is a $S$-free resolution 
of ${\cal I}_{t+1}$, a contradiction.

\item [case $(\rank {\cal I}_{t+1}, \rank {\cal I}_{n-1}) 
= (\binom{n-1}{t}, n-2)$:]
In this case we have $\rank E_{t+1}/{\cal I}_{t+1} =0$.
Since $E_{t+1}/{\cal I}_{t+1}\subset I(c)$ and $I(c)$ is 
torsion-free, we must have $E_{t+1}/{\cal I}_{t+1} =0$. 
Thus $\Ker\phi = E_{t+1}\dirsum {\cal I}_{n-1}$ and 
the Bourbaki sequence (\ref{bourbakiseq0}) is 
\begin{eqnarray*}
     0\To K_n\dirsum U_n \namedTo{\partial_n\dirsum d_n}
                \cdots   \namedTo{\partial_{t+2}\dirsum d_{t+2}}
           & K_{t+1}\dirsum U_{t+1} &  \\
 \namedTo{\partial_{t+1}\dirsum d_{t+1}}
          &E_{t+1}\dirsum E_{n-1}& \namedTo{\phi} I(c)\To 0.\\
\end{eqnarray*}
But then we must have $t \geq n-2$, which contradicts to the 
assumption that $t\leq n-4$.
\end{description}
Now we show $(ii)$ to $(i)$. 
Assume that  (\ref{bourbakiseq0}) is of trivial type. Then
we must have $\beta(G) = p_1(\beta(G))\dirsum p_2(\beta(G))$. From
this we immediately obtain the required result.
\end{proof}

From Lemma~\ref{chainmapBETA}, we immediately have
\begin{Theorem}
\label{theorem:non-trivialcond}
Under the situation of Theorem~\ref{theorem:construction},
the long Bourbaki sequence is of non-trivial type if and only if
\begin{enumerate}
\item [$(i)$] the condition Theorem~\ref{theorem:construction} $(ii)$  holds,
              and
\item [$(ii)$] the submodule $N:= \langle\beta_1,\ldots, \beta_q\rangle$
              of $K_{t+1}\dirsum K_{n-1}(d)$ cannot be decomposed in the form
              of 
	      $N = A\dirsum B$ for some  $(0\ne) A\subset K_{t+1}$ and 
              $(0\ne )B \subset K_{n-1}(d)$
\end{enumerate}
\end{Theorem}

\section{Numerical Characterizations}

Existence of long Bourbaki sequence as in
Theorem~\ref{theorem:construction} and~\ref{theorem:construction2}
only implies that $I$ is a single spot ideal of codimension at most
$3$.  To assure that the codimension is exactly $3$, we need
additional condition. In this section, we give a numerical condition
to assure $\codim I = 3$ for long Bourbaki sequences with approximation
modules $E_{t+1}\dirsum E_{n-1}(d)$.

We assume $n\geq 4$ and $t\leq n-4$, and let $I\subset S$ be a
graded ideal fitting into a long  Bourbaki sequence
\begin{equation}
\label{bourbakiseq}
    0 \To F\namedTo{f} G \namedTo{g} M \namedTo{\phi} I(c) \To 0
\end{equation}
with $M = E_{t+1}\dirsum E_{n-1}(d)$, $F = \Dirsum_{i=1}^p S(-a_i)$ and 
$G=\Dirsum_{i=1}^q S(-b_i)$ are $S$-free modules and $c\in\ZZ$.
As in Theorem~\ref{theorem:construction} we have 
\begin{equation}
\label{rankcondition}
q  = p  + \binom{n-1}{t} + n - 2.
\end{equation}

Now from the sequence (\ref{bourbakiseq}),  we construct 
the mapping cone $C(\alpha, \beta)$ as in (\ref{mappingcone}).
The cone gives a $S$-free resolution $F_\bullet$ of the residue ring $S/I$.
\[
F_\bullet: 0\to F_{n-t}\to\cdots\to F_1 \to F_0 \to S/I \to 0
\]
where 
\begin{eqnarray*}
F_0 & = & S \\
F_1 & = & K_{t+1}(-c)\dirsum K_{n-1}(d-c) = S(-t-1-c)^{\beta_1}\dirsum S(-n+1+d-c)^n  \\
F_2 & = & K_{t+2}(-c)\dirsum K_n(d-c)\dirsum G(-c) \\
    & = & S(-t-2-c)^{\beta_2}\dirsum S(-n+d-c)\dirsum \Dirsum_{i=1}^{q}S(-b_i-c)\\
F_3 & = & K_{t+3}(-c)\dirsum F(-c) = S(-t-3-c)^{\beta_3}\dirsum \Dirsum_{i=1}^{p}S(-a_i-c)\\
F_i & = & K_{t+i}(-c) =S(-t- i -c)^{\beta_{i}} 
\qquad  (4\leq i\leq n-t)\\
\mbox{with} && \beta_i = \binom{n}{t+i}\quad i=1,\ldots, n-t.
\end{eqnarray*}
Notice that this resolution is minmal if and only if matrix
representations of $\alpha$ and $\beta$ only have their entries from
$\mm$.

Now we compute the Hilbert series $Hilb(S/I, \lambda)$ of $S/I$. We have 
\begin{equation}
\label{hilbertSeries}
    Hilb(S/I, \lambda) = \frac{ Q(\lambda)}{(1-\lambda)^n}
\end{equation}
with
\begin{eqnarray*}
\lefteqn{Q(\lambda) =  \sum_{i,j}(-1)^i\beta_{i,j}\lambda^j} \\
           &=& 1 - n \lambda^{n-1+c-d} + \lambda^{n+c-d} 
                + \sum_{i=1}^q\lambda^{b_i+c}
- \sum_{i=1}^p \lambda^{a_i+c}
+ (-1)^t\lambda^c\sum_{i=t+1}^n(-1)^i\binom{n}{i}\lambda^i
\end{eqnarray*}
where $\beta_{i,j}$ are as in $F_i = \Dirsum_{j}S(-j)^{\beta_{ij}}$, $(i=0,\ldots, n-t)$.
(see Lemma~4.1.13 \cite{BH}).

Since we have $\Coh{i}{S/I} = \Coh{i+1}{M}$ for $0\leq i \leq n-4$ by
(\ref{localcohomologies}), we know that $\dim S/I \geq n-3$, i.e.,
$\codim I \leq 3$. To assure that $\codim I \geq 3$ we must have $Q(1)
= Q'(1) = Q''(1) = 0$ (see Corollary~4.1.14(a) \cite{BH}).

\begin{Proposition}
\label{Q(1)=0}
$Q(1) = 0$ holds for all $n,t,c$ and $p$.
\end{Proposition}
\begin{proof}
We compute using (\ref{rankcondition})
\[
 Q(1) = \binom{n-1}{t} - (-1)^{t+1}\sum_{i=t+1}^{n}(-1)^i\binom{n}{i} 
      = \binom{n-1}{t} - \rank E_{t+1} 
      = 0.
\]
where the last equation follows from the Koszul resolution of $E_{t+1}$:
\[
   0\To K_n\To K_{n-1}\To\cdots\To K_{t+1}\To E_{t+1}\To 0\quad\mbox{(exact)}.
\]
\end{proof}

\begin{Proposition}
\label{Q'(1)=0}
$Q'(1) = 0$ holds if and only if 
\[
   \sum_{i=1}^q b_i  - \sum_{i=1}^p a_i 
= n^2 -(2+d)n + c + d + \binom{n-2}{t-1}
                         + \binom{n-1}{t}t.
\]
\end{Proposition}
\begin{proof}
We compute using (\ref{rankcondition})
\[
Q'(1) = -n^2 + (2+d)n - c - d + \sum_{i=1}^q b_i - \sum_{i=1}^p a_i 
           -(-1)^{t+1}\sum_{i=t+1}^n(-1)^i\binom{n}{i}i.
\]
and we know 
\begin{equation}
\label{usefulequation}
(-1)^{t+1}\sum_{i=t+1}^n(-1)^i\binom{n}{i}i
=    (-1)^t\sum_{i=0}^t(-1)^i\binom{n}{i}i
=  \binom{n-2}{t-1} + \binom{n-1}{t}t
\end{equation}
where the last equation is given in Example~2.3 \cite{HT2}.
\end{proof}

Now before we go further, we need to show a combinatorial equation.
\begin{Lemma}
\label{closedF}
\[
(-1)^{t+1}\sum_{i=t+1}^n(-1)^i\binom{n}{i}i^2 
= \binom{n-1}{t}(t+1)^2 - \binom{n-2}{t}(2t+1) - 2\binom{n-3}{t-1}
\]
\end{Lemma}
\begin{proof}
 Let $A = (-1)^{t+1}\sum_{i=t+1}^n(-1)^i\binom{n}{i}i^2$.
A straightforward computation using the binomial coefficient theorem and (\ref{usefulequation}) shows that 
\begin{eqnarray*}
0& = & \sum_{i=2}^{n}(-1)^i\binom{n}{i}i(i-1)  \\
& = & (-1)^{t+1}A 
    + \sum_{i=2}^{t}(-1)^i\binom{n}{i}i(i-1)
    -(-1)^{t+1}\left[ \binom{n-2}{t-1} + \binom{n-1}{t}t\right]
\end{eqnarray*}
Thus we have
\begin{equation}
\label{intermedAns}
A =  \binom{n-2}{t-1} + \binom{n-1}{t}t
    + (-1)^{t}\sum_{i=2}^{t}(-1)^i\binom{n}{i}i(i-1)
\end{equation}
Now we will compute the last term.
Set
\[
\alpha(n,t) = (-1)^{t}\sum_{i = 2}^{t}(-1)^i\binom{n}{i}i(i-1).
\]
and we compute
\begin{eqnarray*}
\lefteqn{\alpha(n,t) + \alpha(n-1,t-1)}\\
&= & 2t^2\binom{n-1}{t}  -
2(-1)^{t}\sum_{i=2}^{t}(-1)^i\binom{n-1}{i}i
     + 2(-1)^{t}(n-1) - \alpha(n-1, t).
\end{eqnarray*}
Also we have
\[
	\alpha(n-1, t-1) + \alpha(n-1, t)
	= t(t-1)\binom{n-1}{t}
\]
Thus by using (\ref{usefulequation})
\begin{eqnarray*}
\alpha(n,t)
& = & t(t+1)\binom{n-1}{t} -
2(-1)^{t}\sum_{i=2}^{t}(-1)^i\binom{n-1}{i}i
	+  2(-1)^{t}(n-1) \\
& = & t(t+1)\binom{n-1}{t} - 2
	\left[
		t\binom{n-2}{t} + \binom{n-3}{t-1}
              + (-1)^{t}(n-1)
	\right]\\
&   &
	+  2(-1)^{t}(n-1) \\
& = & t(t+1)\binom{n-1}{t} - 2t\binom{n-2}{t} -
2\binom{n-3}{t-1}\\
\end{eqnarray*}
Substiting this into  (\ref{intermedAns}), we obtain the desired result.
\end{proof}

\begin{Proposition}
\label{Q''(1)=0}
$Q''(1)=0$ holds if and only if 
\begin{eqnarray*}
\sum_{i=1}^q b_i^2  - \sum_{i=1}^p a_i^2
&=& n^3 - (3+2d)n^2  + (d^2+4d+1)n - c^2 - d^2  \\
& &
 + \binom{n-1}{t}(t+1)^2 
    - \binom{n-2}{t}(2t+1) - 2\binom{n-3}{t-1}
\end{eqnarray*}
\end{Proposition}
\begin{proof}
We compute 
\begin{eqnarray*}
Q''(1) &=& c + d -c^2 + d^2 - 2cd + (4c-3-d^2-5d+2cd)n + (4-2c+2d)n^2 - n^3 \\
       & & -(2c-1)\left[\binom{n-2}{t-1} + \binom{n-1}{t}t\right] 
             + \sum_{i=1}^{q}b_i^2 - \sum_{i=1}^p a_i^2\\
       & &   + (2c -1)\left(\sum_{i=1}^{q}b_i - \sum_{i=1}^{p}a_i\right) 
             + (-1)^t\sum_{i=t+1}^{n}(-1)^i\binom{n}{i}i^2 \\
       &=&  c^2 +d^2 -(d^2+4d+1)n + (3+2d)n^2 - n^3  \\
       & &  + \sum_{i=1}^q b_i^2  - \sum_{i=1}^p a_i^2
            + (-1)^t\sum_{i=t+1}^n(-1)^i\binom{n}{i}i^2\\ 
\end{eqnarray*}
where the last equation is by  Proposition~\ref{Q'(1)=0}. 
Then by Lemma~\ref{closedF} we obtain the desired result.
\end{proof}

To summerize, we obtain 
\begin{Theorem}
\label{theorem:numerical}
Let $n\geq 4$ and $t\leq n-4$. Assume that we have the following long Bourbaki sequence
\[
0\To \Dirsum_{i=1}^p S(-a_i)\To \Dirsum_{i=1}^q S(-b_i)\To E_{t+1}\dirsum E_{n-1}(d)\To I(c)\To 0
\]
with $I \subset S$ a graded ideal and $c\in\ZZ$. Then we have $\codim
I \leq 3$ and the equality holds if and only if
\begin{enumerate}
\item $q = p + \binom{n-1}{t} + n - 2$;
\item $\displaystyle{\sum_{i=1}^q }b_i - \displaystyle{\sum_{i=1}^p} a_i
 = n^2 - (2+d)n + c + d + \binom{n-2}{t-1} + \binom{n-1}{t}t$;  
\item and \begin{eqnarray*}
\displaystyle{\sum_{i=1}^q} b^2_i - \displaystyle{\sum_{i=1}^p} a^2_i
               &=& n^3 - (3+2d)n^2 + (d^2+4d+1)n - c^2 - d^2 \\
               & & + \binom{n-1}{t}(t+1)^2 -\binom{n-2}{t}(2t+1) - 2\binom{n-3}{t-1}
      \end{eqnarray*}
\end{enumerate}
\end{Theorem}

%

\section{Examples}

\begin{Example}
\label{example1}
We first give an application of Theorem~\ref{theorem:construction2}. Namely,
a single spot ideal $I$ with approximation module $E_{t+1}$. 
Let $t = 1$ and $n =6$. Then ${\cal A} = \{A_1,A_2,A_3,A_4,A_5,A_6\}
 = \partial_2^*(E_2^*) \subset K_2^*$ where 
\begin{eqnarray*}
A_1 & = &  x_1 e^*_{16} + x_2 e^*_{26} + x_3 e^*_{36} + x_4 e^*_{46} + x_5 e^*_{56}
\\
A_2 & = & -x_1 e^*_{15} - x_2 e^*_{25} - x_3 e^*_{35} - x_4 e^*_{45} + x_6 e^*_{56}
\\
A_3 & = &  x_1 e^*_{14} + x_2 e^*_{24} + x_3 e^*_{34} - x_5 e^*_{45} - x_6 e^*_{46}  
\\
A_4 & = & -x_1 e^*_{13} - x_2 e^*_{23} + x_4 e^*_{34} + x_5 e^*_{35} + x_6 e^*_{36}  
\\
A_5 & = &  x_1 e^*_{12} - x_3 e^*_{23} - x_4 e^*_{24} - x_5 e^*_{25} - x_6 e^*_{26}  
\\
A_6 & = &  x_2 e^*_{12} + x_3 e^*_{13} + x_4 e^*_{14} + x_5 e^*_{15} + x_6 e^*_{16}.
\end{eqnarray*}
Now let $a \in {\cal A}$ and $\beta_i$ ($i=1,\ldots, 6$) be as follows:
\begin{eqnarray*}
a & = & x_6 A_1 -x_5 A_2 + x_4 A_3 \\
  & = & x_1x_4 e^*_{14} + x_1x_5 e^*_{15} +x_1x_6 e^*_{16} + x_2x_4 e^*_{24}
       x_2x_5 e^*_{25} + x_2x_6 e^*_{26} +x_3x_4 e^*_{34} + x_3x_5 e^*_{35}\\
  &   & + x_3x_6 e^*_{36} \\
\beta_1 & = & e_{12},\quad \beta_2 = e_{13},\quad \beta_3 = e_{23},
                     \quad \beta_4 = e_{45},\quad \beta_5 = e_{46},
                     \quad \beta_6 = e_{56}
\end{eqnarray*}
Then we have
\[
 \Ker(a : K_2\dirsum K_5 \to S)
  = \langle\beta_1,\ldots,\beta_6\rangle + E_3
\]
and, for the map $\beta : \Dirsum_{i=1}^6 S(-2) \to K_2\dirsum K_5$ such that
$\beta(m_i) = \beta_i$ $(i=1,\ldots, 6)$ where $\{m_i\}$ is a free basis, 
we obtain the diagram
\[
\begin{CD}
   @.    0      @.   0     \\
@.       @VVV        @VVV  \\
0 @>>> \left\langle
\begin{array}{c}
          x_3m_1-x_2m_2+x_1m_3,\\
          x_6m_4-x_5m_5+x_4m_6 
         \end{array}
         \right\rangle
                @>>> \langle m_1,\ldots,m_6\rangle \\
@.       @V{Res(\beta)}VV
                      @V{\beta}VV \\
0 @>>> \langle\beta_1,\ldots,\beta_6\rangle\cap E_3
                 @>>> \langle\beta_1,\ldots,\beta_6\rangle \\
@.       @VVV          @VVV  \\
    @.   0       @.     0    \\
\end{CD}
\]
and $a$ defines a degree $0$ homomorphism from $E_2$ to $S$.
Then we obtain the long Bourbaki sequence
\[
   0 \To S^2(-3) \namedTo{f} S^6(-2) \namedTo{g} E_2 \namedTo{\varphi} I \To 0
\]
where 
\[
\begin{array}{cccc}
f : & S^2(-3) = S n_1 \dirsum S n_2 & \To &  S^6(-2) = S m_1 \dirsum\cdots\dirsum S m_6 \\
    &           n_1                 &\longmapsto  &  x_3 m_1 - x_2 m_2 + x_1 m_3 \\
    &           n_2                 &\longmapsto  &  x_6 m_4 - x_5 m_5 + x_4 m_6 \\
    &                               &             &                               \\
g : &  S^6(-2) = S m_1 \dirsum\cdots\dirsum S m_6 & \To  & E_2 \\
    &           m_i                 &\longmapsto  &  \partial_2(\beta_i) \qquad (i=1,\ldots,6) \\
\mbox{and} &                        &             &                               \\
\varphi &   E_2               & \To         & I  \\
        & \partial_2(e_{ij})  & \longmapsto & x_ix_j \\
        &                     &\mbox{for }& (i,j) = (1,2),(1,3),(2,3),(4,5),(4,6),(5,6)\\
        & \partial_2(e_{ij})  & \longmapsto & 0 \\
        &                     &\mbox{for }& (i,j) \ne (1,2),(1,3),(2,3),(4,5),(4,6),(5,6)\\
\end{array}
\]
and we obtain 
$I = (x_1x_4,x_1x_5, x_1x_6,x_2x_4,x_2x_5,x_2x_6,x_3x_4, x_3x_5, x_3x_6)
       = (x_1,x_2,x_3)(x_4,x_5,x_6)$, a codimension~3 
single spot ideal of type $(1,K)$.
\end{Example}

\begin{Example}
\label{example2}
We continue to consider the situation in Example~\ref{example1}. 
As an application of Theorem~\ref{theorem:construction}, we can see
that the same ideal fits into a long Bourbaki sequence with approximation
module $E_{t+1}\dirsum E_{n-1} = E_1 \dirsum E_5$.
In this case, we must also  consider 
${\cal B} = \{B_{ij} \mid 1\leq i <j \leq 6 \}
 = \partial_5^*(E_5^*) \subset K_5^*$ 
where $B_{ij} = (-1)^{i}x_je_{[6]\backslash i}^*
              - (-1)^{j}x_ie_{[6]\backslash j}^*$.
Then we set $a \in {\cal A}$ as in Example~\ref{example1} and 
\[
b = -x_1^2x_2x_4 B_{14} =x_1^2x_2x_4^2 e^*_{23456} + x_1^3x_2x_4 e^*_{12356}
\in {\cal B}.
\]
Also we set $\beta_1,\ldots, \beta_6$ to be the same as 
those in Example~\ref{example1} and 
$\beta_7 =x_1x_2x_4 e_{14} - e_{23456}$,
$\beta_8 =x_1^2x_2 e_{14} - e_{12356}$,
$\beta_9    = e_{13456}$,
$\beta_{10} = e_{12456}$,
$\beta_{11} = e_{12346}$,
$\beta_{12} = e_{12345}$.
Notice that $\{\beta_i\}_{i=1}^{12}$ satisfies the condition of
Theorem~\ref{theorem:non-trivialcond}.
Then $\varphi = (a,b)$ defines a degree $0$ map on $E_2\dirsum E_5$, and we have
\[
  \Ker\varphi = \langle\beta_1,\ldots,\beta_{12}\rangle + E_3\dirsum E_6
\]
and the diagram
\[
\begin{CD}
   @.    0     @.     0  \\
@.       @VVV        @VVV \\
0 @>>> \left\langle
          \begin{array}{l}
            x_3m_1-x_2m_2+x_1m_3,\;
            x_6m_4-x_5m_5+x_4m_6,\\
          x_1m_7-x_4m_8+x_2m_9-x_3m_{10}-x_5m_{11}+x_6m_{12}
          \end{array}
         \right\rangle
              @>>> \langle m_1,\ldots,m_{12}\rangle \\
@.       @V{Res(\beta)}VV    
                   @V{\beta}VV \\
0 @>>>  \langle\beta_1,\ldots,\beta_6\rangle
          \cap E_3\dirsum E_6
                   @>>> \langle\beta_1,\ldots\beta_{12}\rangle\\
@.        @VVV           @VVV   \\
    @.    0        @.     0
\end{CD}
\]
Then we have a long Bourbaki sequence of non-trivial type
\[
   0 \To S^2(-3)\dirsum S(-6) 
    \namedTo{f} 
         S^6(-2)\dirsum S^6(-5)
    \namedTo{g} 
         E_2\dirsum E_5
    \namedTo{\varphi} I \To 0
\]
where 
\[
\begin{array}{cccc}
f : & S^2(-3)\dirsum S(-6) = S n_1 \dirsum S n_2 \dirsum S n_3 & \To &  S^6(-2) = S m_1 \dirsum\cdots\dirsum S m_6 \\
    &           n_1                 &\longmapsto  &  x_3 m_1 - x_2 m_2 + x_1 m_3 \\
    &           n_2                 &\longmapsto  &  x_6 m_4 - x_5 m_5 + x_4 m_6 \\
    &           n_3                 &\longmapsto  &  x_1 m_7 - x_4 m_8 + x_2 m_9 \\
    &                               &             & - x_3 m_{10} -x_5 m_{11} + x_6 m_{12} \\
    &                               &             &                               \\
g : &  S^6(-2) = S m_1 \dirsum\cdots\dirsum S m_6 & \To  & E_2\dirsum E_5 \\
    &           m_i                 &\longmapsto  &  \bar\partial(\beta_i) \qquad (i=1,\ldots,12) \\
    &                               &             & \mbox{where }\bar\partial = \partial_2\dirsum\partial_5\\
\mbox{and} &                        &             &                               \\
\varphi &   E_2               & \To         & I  \\
        & \partial_2(e_{ij})  & \longmapsto & x_ix_j \\
        &                     &\mbox{for }& (i,j) = (1,2),(1,3),(2,3),(4,5),(4,6),(5,6)\\
        & \partial_2(e_{ij})  & \longmapsto & 0 \\
        &                     &\mbox{for }& (i,j) \ne (1,2),(1,3),(2,3),(4,5),(4,6),(5,6)\\
        & \partial_5(e_{23456}) & \longmapsto & x_1^2x_2x_4^2 \\
        & \partial_5(e_{12356}) & \longmapsto & x_1^3x_2x_4 \\
        & \partial_5(e_{ijklm}) & \longmapsto & 0\quad \mbox{otherwise} \\
\end{array}
\]
and the ideal $I$ is the same as that in Example~\ref{example1}.
We can also check that this sequence satisfies the numerical condition in Theorem~\ref{theorem:numerical}
\end{Example}

\begin{Example}
By Corollary~\ref{cor:depthzero}, we do not have a long Bourbaki
sequence with an approximation module $E_1$ and a codimension~$3$
generalized CM ideal $I$. However, there are long Bourbaki sequences
with approximation modules $E_1\dirsum E_5(d)$ for $d\in\ZZ$, which 
is an application of Theorem~\ref{theorem:construction}.
Let $k=1$ and $n=6$. Then 
\begin{eqnarray*}
{\cal A} & =&  \langle
                   x_1e_1^* + \cdots + x_6 e_6^* 
               \rangle\\
{\cal B} &= & \langle B_{ij} = (-1)^ix_je_{[6]\backslash i}^*
- (-1)^jx_ie_{[6]\backslash j}^*  :
		  1\leq i < j \leq 6 
\rangle.
\end{eqnarray*}
Let $\varphi = (a,b)$ be
\begin{eqnarray*}
a & = & x_1^3e_1^* + x_1^2x_2 e_2^* + x_1^2x_3 e_3^* + 
        x_1^2x_4 e_4^* + x_1^2x_5 e_5^* + x_1^2x_6 e_6^* \\
b & = & -x_2^5 B_{56} + x_6^5 B_{23}\\
  & = & x_2^5x_6e^*_{12346} + x_2^5x_5 e^*_{12345}
         +x_3x_6^5 e^*_{13456} + x_2x_6^5 e^*_{12456} \\
\end{eqnarray*}
and set $\beta_i\in K_1\dirsum K_5(1)$ to be as follows:
\begin{eqnarray*}
\beta_1 & = & - x_6e_{12345}  + x_5 e_{12346} \\
\beta_2 & = &  x_6^5 e_3 - x_1^2 e_{13456} \\
\beta_3 & = &  x_6^5 e_2 - x_1^2 e_{12456} \\
\beta_4 & = &  x_2^4x_5 e_2 - x_1 e_{12345} \\
\beta_5 & = &  x_2^4x_6 e_2 - x_1^2 e_{12346} \\
\beta_6 & = & -x_6^4 e_{12346} + x_2^4 e_{12456} \\
\beta_7 & = & e_{23456} \\
\beta_8 & = & e_{12356}.
\end{eqnarray*}
Notice that $\beta_i \notin E_2\dirsum E_6(1)$ for all $i$, i.e., the
condition in Theorem~\ref{theorem:non-trivialcond} is satisfied.
Then we can check 
\begin{enumerate}
\item $\Ker (\varphi: K_1\dirsum K_5(1)\to S(-6) 
      = \langle \beta_1,\ldots, \beta_8\rangle + E_2\dirsum E_6(1)$ and 
     $\varphi$ is a degree $8$ homomorphism, and 
\item the diagram
\[
\begin{CD}
   @.    0                @.    0                  \\
@.       @VVV             @VVV                     \\
0  @>>>  \Ker Res(\beta)  @>>>  \Ker\beta   @>>> 0 \\
@.       @VVV                   @VVV           @.  \\
0  @>>>  F'               @>>>  G           @.     \\
@.       @V{Res(\beta)}VV       @V{\beta}VV    @.  \\
0 @>>>  \langle\beta_1,\ldots,\beta_8\rangle\cap E_2\dirsum E_6(1)
                          @>>>  \langle\beta_1,\ldots,\beta_8\rangle \\
@.        @VVV                  @VVV           @.   \\
    @.    0               @.    0           @.      \\
\end{CD}
\]
where 
\begin{eqnarray*}
G & = & \langle m_1,\ldots, m_8\rangle
=S(-5)\dirsum S^4(-6)\dirsum S(-8)\dirsum S^2(-4) \\
F' & = & \left\langle 
\begin{array}{l}
-x_1^2 m_1 + x_6 m_4 - x_5 m_5,
                 x_2^4 m_3 - x_6^4 m_5 + x_1^2 m_6,\\
                -x_1^2m_1 -x_2 m_2 + x_3 m_3 - x_1^3 m_7 + x_1^2x_4 m_8
             \end{array}
         \right\rangle \\
\Ker\beta & = & \Ker Res(\beta) = \langle -x_1^2 m_1 + x_6 m_4 - x_5 m_5,
                 x_2^4 m_3 - x_6^4 m_5 + x_1^2 m_6
         \rangle \\
\end{eqnarray*}
\end{enumerate}
Thus we have a long Bourbaki sequence
\[
  0 \to F \namedTo{f} G \namedTo{g} E_1\dirsum E_5(1) \namedTo{\phi} I(2) \to 0
\]
where $g(m_i) = \beta_i$, $i=1,\ldots, 8$, and $F = S(-10)\dirsum S^2(-7)
= \langle u, v, w\rangle$ with 
$f(u) = x_2^4 m_3 - x_6^4 m_4 + x_1^2 m_6$,
$f(v) = -x_1^2 m_1 + x_6 m_4 - x_5m_5$
and 
$f(w) = -x_1^2 m_1 - x_2 m_2
+ x_3 m_3 - x_1^3 m_7 + x_1^2x_4 m_8$. The map $\phi$ is as follows:
$\phi(x_i) = x_ix_1^2$ $(i=1,\ldots, 6)$, $\phi(\partial_5(e_{12345})) = x_2^5x_5$,
$\phi(\partial_5(e_{12346})) = x_2^5x_6$, $\phi(\partial_5(e_{12356})) = 0$,
$\phi(\partial_5(e_{12456})) = x_2x_6^5$, $\phi(\partial_5(e_{13456})) 
   = x_3x_6^5$, and $\phi(\partial_5(e_{23456})) = 0$. 
The ideal is $I = \Im\varphi = x_1^2\mm + (x_2^5x_6, x_2^5x_5, x_3x_6^5, x_2x_6^5)$.
Finally we can check that this Bourbaki sequence satisfies the numerical
condition of Theorem~\ref{theorem:numerical}.
\end{Example}



\begin{thebibliography}{99}

\bibitem{Ama1} M.\ Amasaki, Basic sequences of homogeneous ideals in 
polynomial rings, J.\ Algebra ~190, pp329-360, 1997.


\bibitem{Bour} N.\ Bourbaki, {\em Elements of Mathematics, Commutative Algebra, Chapter 1-7},
Springer, 1989.

\bibitem{BH} W.\ Bruns and   J.\ Herzog, {\em Cohen-Macaulay Rings}, Revised
version, Cambridge University Press, 1998.



\bibitem{HT2} J.\ Herzog and Y.\ Takayama, Approximations of Generalized
Cohen-Macaulay Modules, preprint,  2002.

\bibitem{Kunz} E.\ Kunz, {\em Introduction to Commutative Algebra and 
Algebraic Geometry}, Birkh\"auser, 1985.


 




\bibitem{F} H.\ Flenner, Die S\"atze von Bertini f\"ur lokale Ringe, Math.
Ann. {\bf 299}, 97--111, 1977.



\end{thebibliography}
\end{document}